\documentclass[12pt,]{article}
\usepackage{authblk}
\usepackage{natbib}
\usepackage{amssymb,amsmath}
\usepackage[utf8x]{inputenc}
\usepackage{authblk}
\providecommand{\authorinfo}[1]{\noindent{#1}\vspace{2em}}
\providecommand{\keywords}[1]{\textbf{Keywords: } #1}
\providecommand{\suppmaterial}[1]{\noindent\textbf{Supplementary material: } #1\vspace{1em}}
\providecommand{\abstractm}[1]{\noindent\textbf{Abstract: } #1\vspace{1em}}
\usepackage{array}
\usepackage{tabularx}
\usepackage{multirow}
\usepackage{quoting}
\usepackage{setspace}
\usepackage[unicode=true]{hyperref}
\usepackage{color}
  \definecolor{mygray}{rgb}{0.9,0.9,0.9}
\usepackage{textcomp}
\usepackage{listings}
\usepackage{hyperref}
\usepackage{graphicx}
\graphicspath{{}}
\usepackage[normal]{subfigure}

\usepackage{array} 
\usepackage{tabularx} 
\usepackage{multirow} 
\usepackage{rotating} 
\usepackage{lscape} 
\title{A generalized planar conchoid}
\author[1]{Ludger O. Suarez-Burgoa}
\affil[1]{Universidad Nacional de Colombia}
\date{\today}

\begin{document}
\maketitle

\authorinfo{%
Universidad Nacional de Colombia, Faculty of Mines\\
ORCID ID: FA, 0000-0002-9760-0277\\
Present address: FA, Department of Civil Engineering, Cl. 65 \# 78-28, Medell\'in, AN 050034, Colombia\\
Corresponding author (e-mail: losuarezb@unal.edu.co)
}

\abstractm{
This note presents the definition of a proposed generalization of the conchoid at the plane. Known conchoids, such as the Nicomedes and the Lima\c{c}on of Pascal are part of this set. Following the definition, one can generate other conchoids. Examples are generated using of a computer code that is available openly for download. In addition, two step-by-step examples are described by detail, the first one which presents the results in calculation tables.
}


\suppmaterial{%
The computation code package, named \textsf{genPlanarConchoid}, is available at \\ \url{https://github.com/losuarezburgoa/genPlanarConchoid}.
}
\newpage


\section{Definition}

\begin{figure}[!ht]
  \centering
  \subfigure[Parts for its definition]{\label{subfig:partsOfGCdefinitionV2} \includegraphics[width=0.45\textwidth]{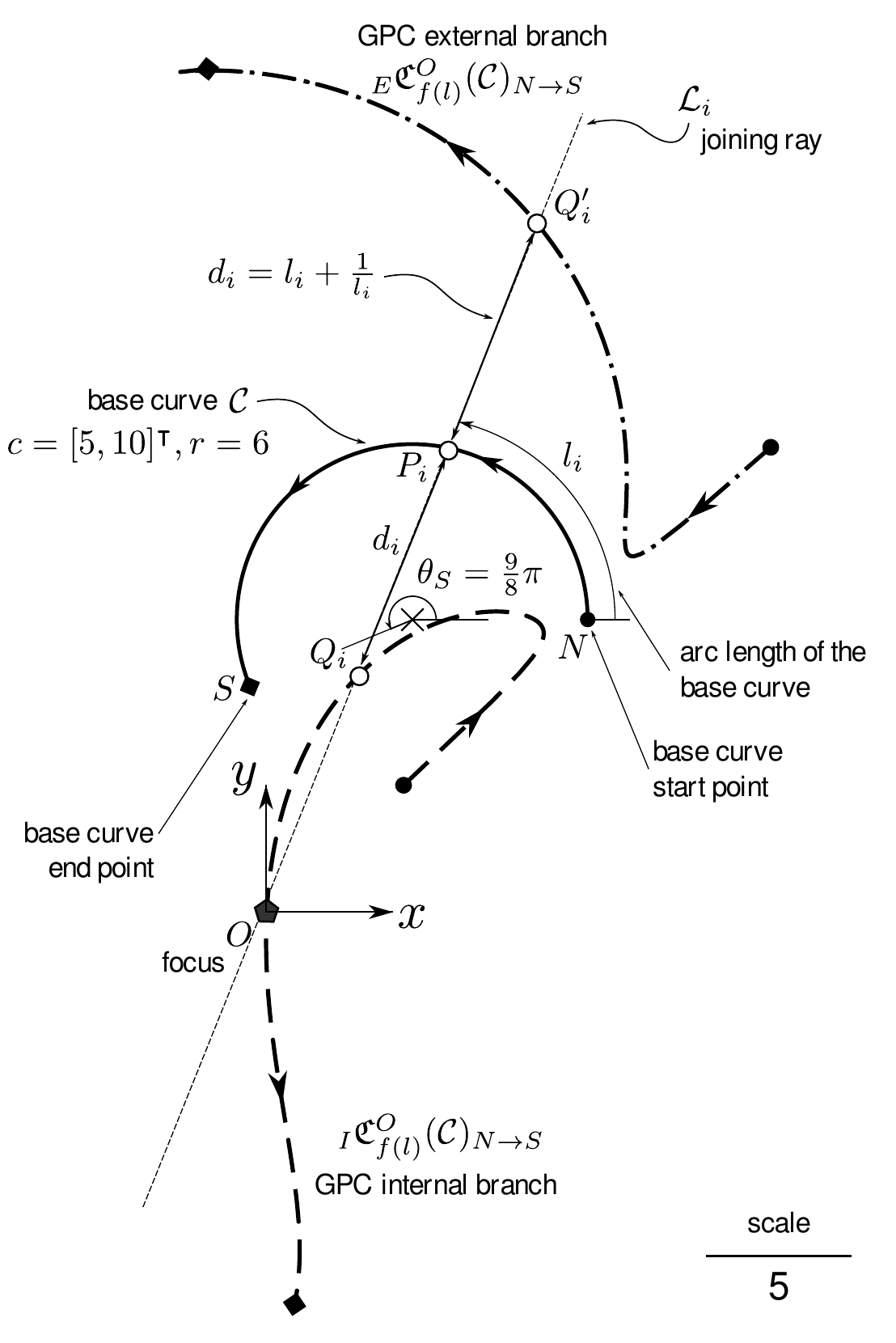}}
  \quad
  \subfigure[The curves of interest]{\label{subfig:resultingGPCofOfDefinition} \includegraphics[width=0.45\textwidth]{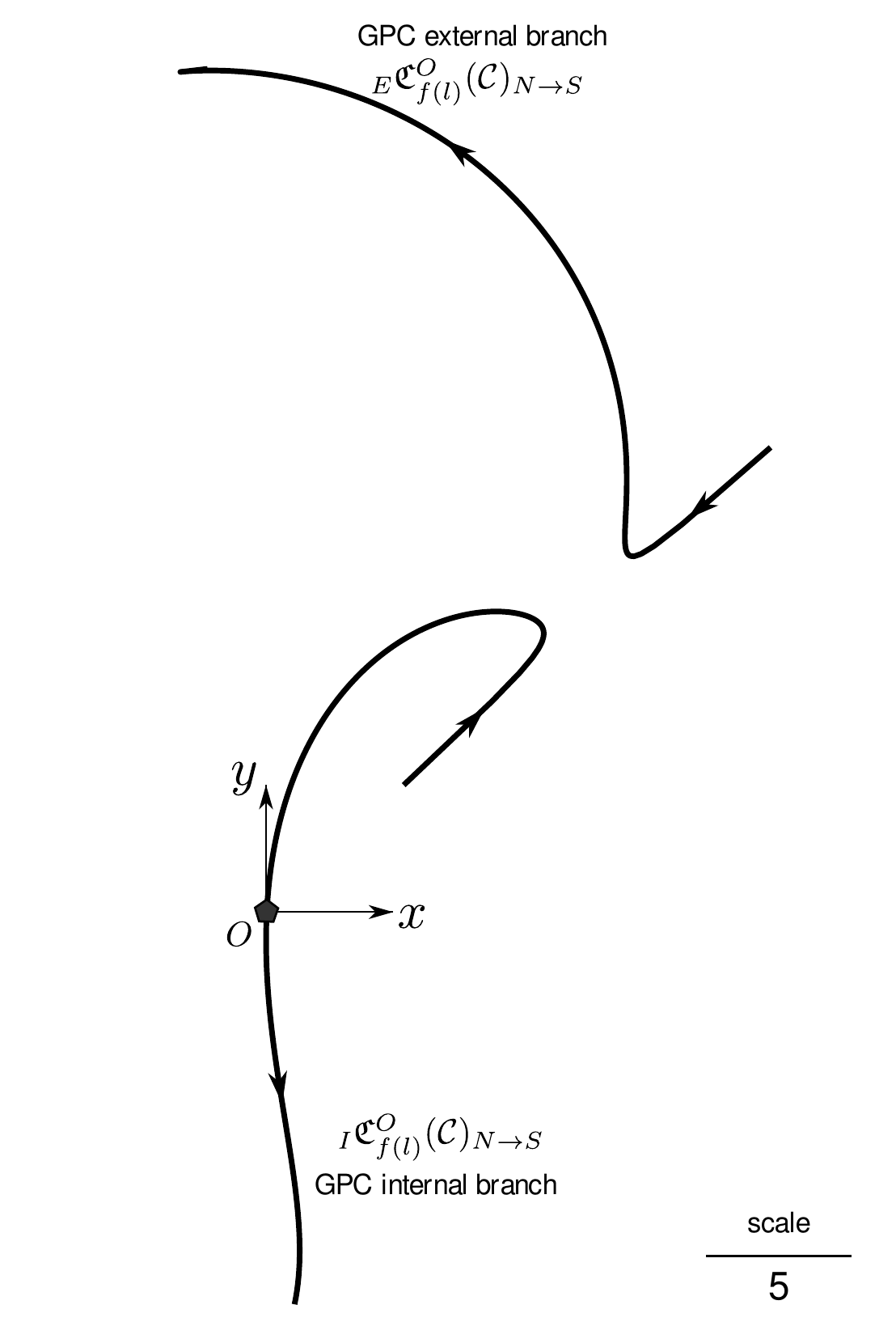}}\\
  \caption{General planar conchoid (GPC); in this case, it is \(\mathfrak{C}_{f(l)=l + \frac{1}{l}}^{O=[0, 0]^\intercal}(\mathcal{C}: c=[5, 10]^\intercal, r=6)_{\theta_N=0 \rightarrow \theta_S = \frac{9}{8}\pi}\) (Cont.).}
  \label{fig:partsOfGCdefinition}
\end{figure}

Let be \(O\) a fixed point called a focus and let \(\mathcal{L}_i\) be a set of lines, where the lines are required to pass through \(O\) and intersect a curve \(\mathcal{C}\) at points \(P_i\). The geometric locus of points \(Q_i\) and \(Q'_i\) on \(\mathcal{L}_i\) such that a variable offset Euclidean distance
\(d_i = \overline{P_iQ_i} = \overline{P_iQ'_i}\) (for \(d_i \in \mathbb{R}^+ > 0\)) responds to a function \(f(l_i)\), which depends on the arc-length (\(l_i\)) measured from a starting point \(N\) and directly passes through any \(P_i\) and ending at point \(S\) in \(\mathcal{C}\), is defined here as a \emph{generalized planar conchoid} (GPC) and is denoted as \(\mathfrak{C}_{f(l)}^{O}(\mathcal{C})\). Because \(N\) and \(S\) are the starting and ending points of \(\mathcal{C}\), \(\mathcal{C}\) is finite and has a direction.
The notation \(\mathfrak{C}_{f(l)}^{O}(\mathcal{C})_{N \rightarrow S}\) is read as:
\begin{quoting}
a generalized planar conchoid at focus \(O\) with base curve \(\mathcal{C}\) from \(N\) to \(S\), and is based on function \(f(l)\).
\end{quoting}
Figure \ref{fig:partsOfGCdefinition} shows the conceptual scheme of GPC with the names of each of the parts.

Point \(O\) is (as mentioned) the \emph{focus}. The curve \(\mathcal{C}\) is the \emph{base curve}. Lines \(\mathcal{L}_i\) are the rays that define two segments from point \(P_i\) to \(Q_i\) and from \(P_i\) to \(Q'_i\), which are called the \emph{interior branch} (\(\overline{P_i Q_i}\)) and the \emph{exterior branch} (\(\overline{P_i Q'_i}\)), respectively. The distance \(d_i\), which is equal to the interior and exterior branch distances, is called the \emph{distance offset}. The distance from \(N\) to \(P_i\) through \(\mathcal{C}\) is the \emph{arc length} (\(l_i\)). The function \(f(l_i)\) is called the \emph{offset function} or \emph{distance function}. Point \(N\) is the starting point of the directed curve \(\mathcal{C}\), upon which the \emph{arc length} is measured, and point \(S\) is the ending point of \(\mathcal{C}\).

The GPC (\(\mathfrak{C}_{f(l)}^{O}(\mathcal{C})_{N \rightarrow S}\)) is composed of the following mathematical objects.
\begin{enumerate}
\item The \emph{focus} \(O\) is a point represented by a column vector with a size of 2 \(\times\) 1.
\item The \emph{base curve} \(\mathcal{C}\) is represented by the following objects.
\begin{enumerate}
\item The function of the curve: \(c(x,y)\).
\item A starting point \(N\), where \(N \in c(x,y)\).
\item An ending point \(S\), where \(S \in c(x,y)\).
\item The function that defines the arc length of the base curve: \(l = g(x,y)\).
\end{enumerate}
\item The function relating the distance \(d\) and the arc length: \(d = f(l)\).
\end{enumerate}

In the next section, some examples that are used create a GPC that follows the defined rules are shown.

\section{How to create a GPC}

The creation of a GPC is straightforward when following the definition described in the previous section, but for a rapid implementation, a computation code in any programming language can be employed.

A particular difficulty may arise in the calculation of the \emph{arc length} of the curve \(\mathcal{C}\). Therefore, for each \(\mathcal{C}\) curve, a function for the arc-length estimation should be defined. For some particular planar curves, the arc length is not exactly defined, for example, when an ellipse is used as a base curve.

In this text, two Octave/MATLAB functions are implemented for the case when the base curves (\(\mathcal{C}\)) are a line and an arc circle; those functions are \textsf{linegenconchoid} and \textsf{circarcgenconchoid}, respectively. Both functions return a data structure that is loaded into a plotting function (called \textsf{plotgenconchoid}) to obtain graphical representations.

\begin{figure}[!ht]
  \centering
  \subfigure[{\(\mathfrak{C}_{f(l)=2}^{O=[0, 0]^\intercal}(\mathcal{C}: y=1)_{N_x=-3  \rightarrow S_x= 3}\)}]{\label{subfig:linConchConstFun} \includegraphics[width=0.45\textwidth]{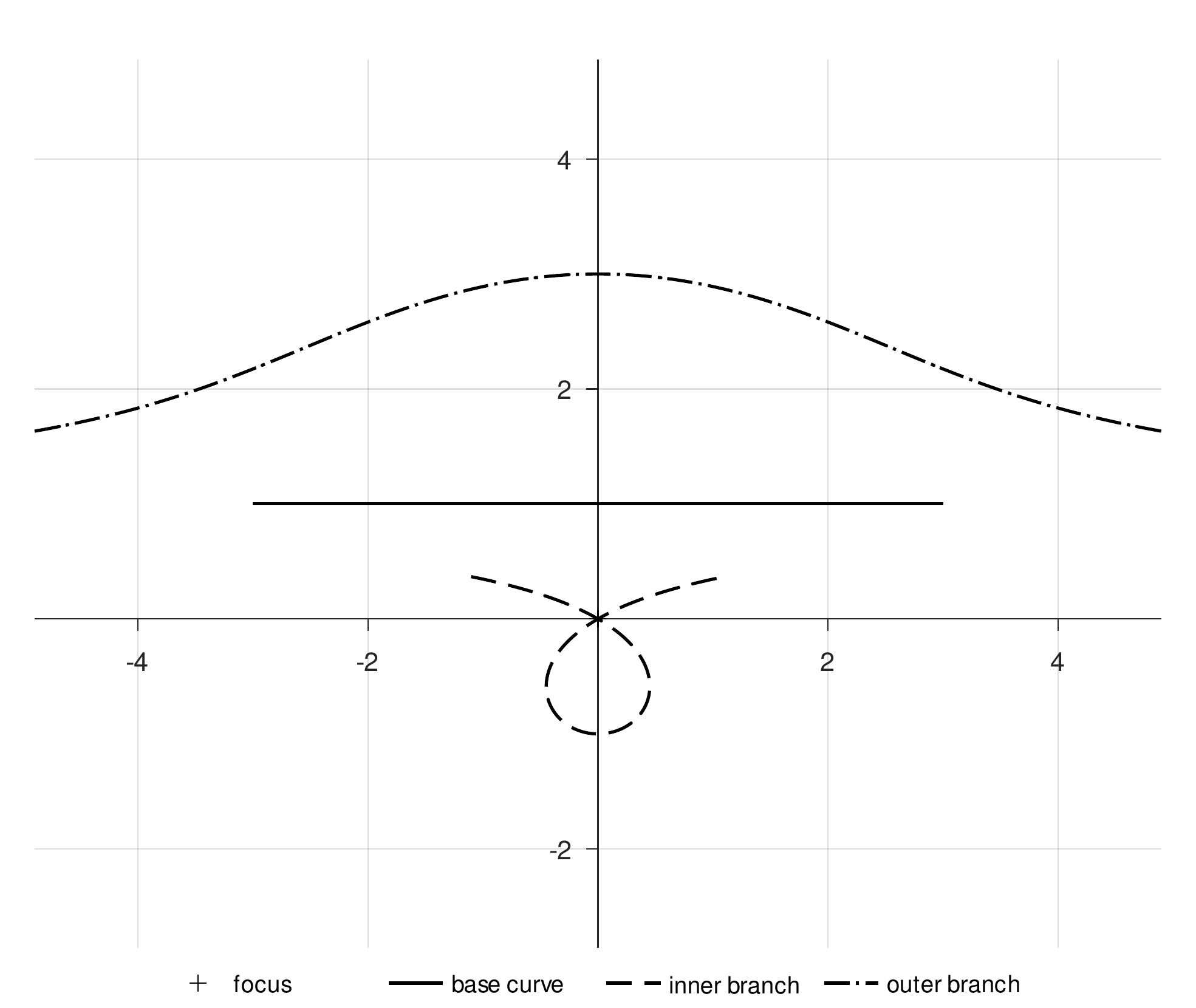}}\quad
  \subfigure[{\(\mathfrak{C}_{f(l)=l}^{O=[0, 0]^\intercal}(\mathcal{C}: y=1)_{N_x=-3  \rightarrow S_x= 3}\)}]{\label{subfig:linConchLinFun} \includegraphics[width=0.45\textwidth]{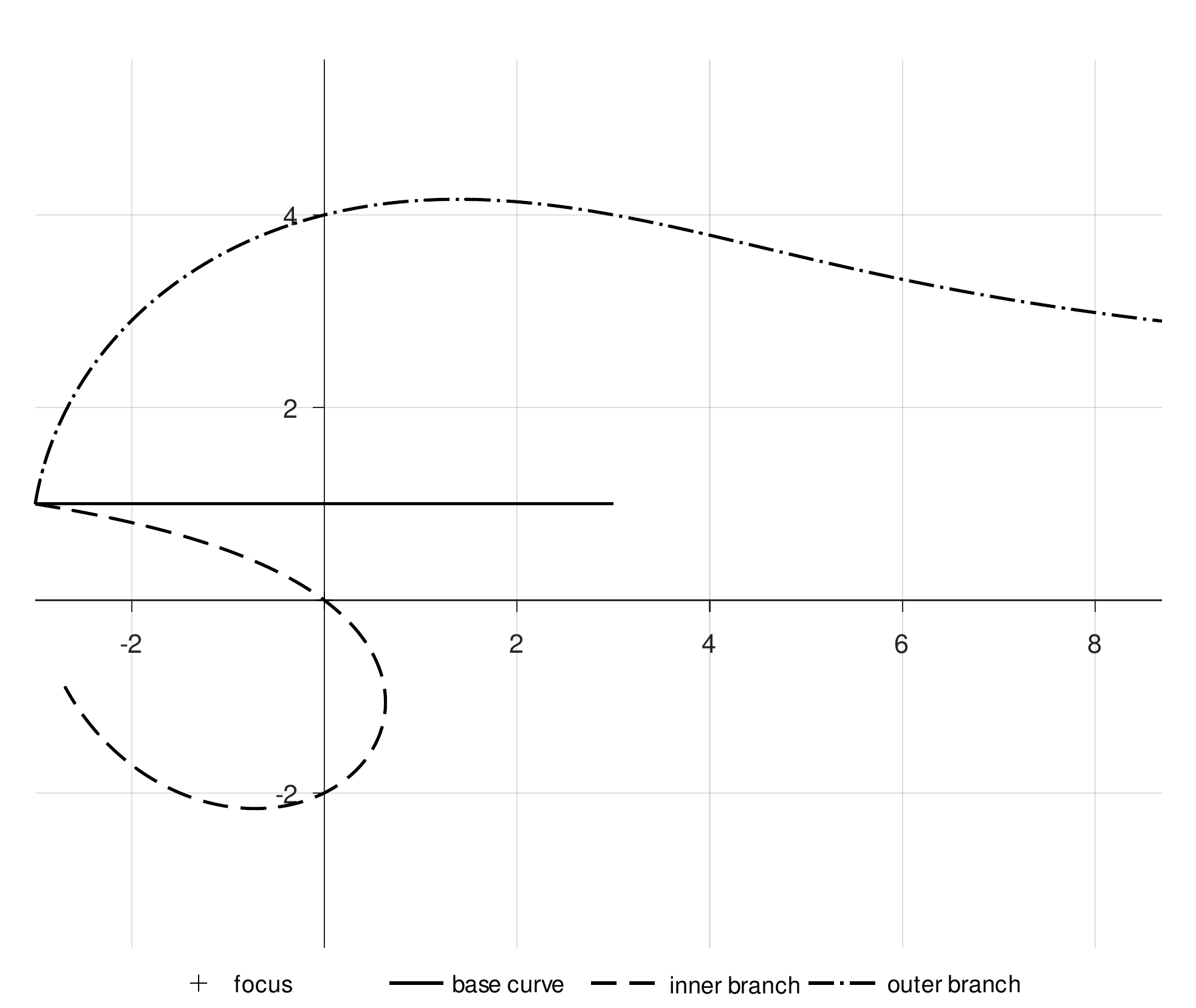}}\\
  \subfigure[{\(\mathfrak{C}_{f(l)=\sin{l}}^{O=[0, 0]^\intercal}(\mathcal{C}: y=1)_{N_x=-4  \rightarrow S_x= 4}\)}]{\label{subfig:linConchSinFun} \includegraphics[width=0.45\textwidth]{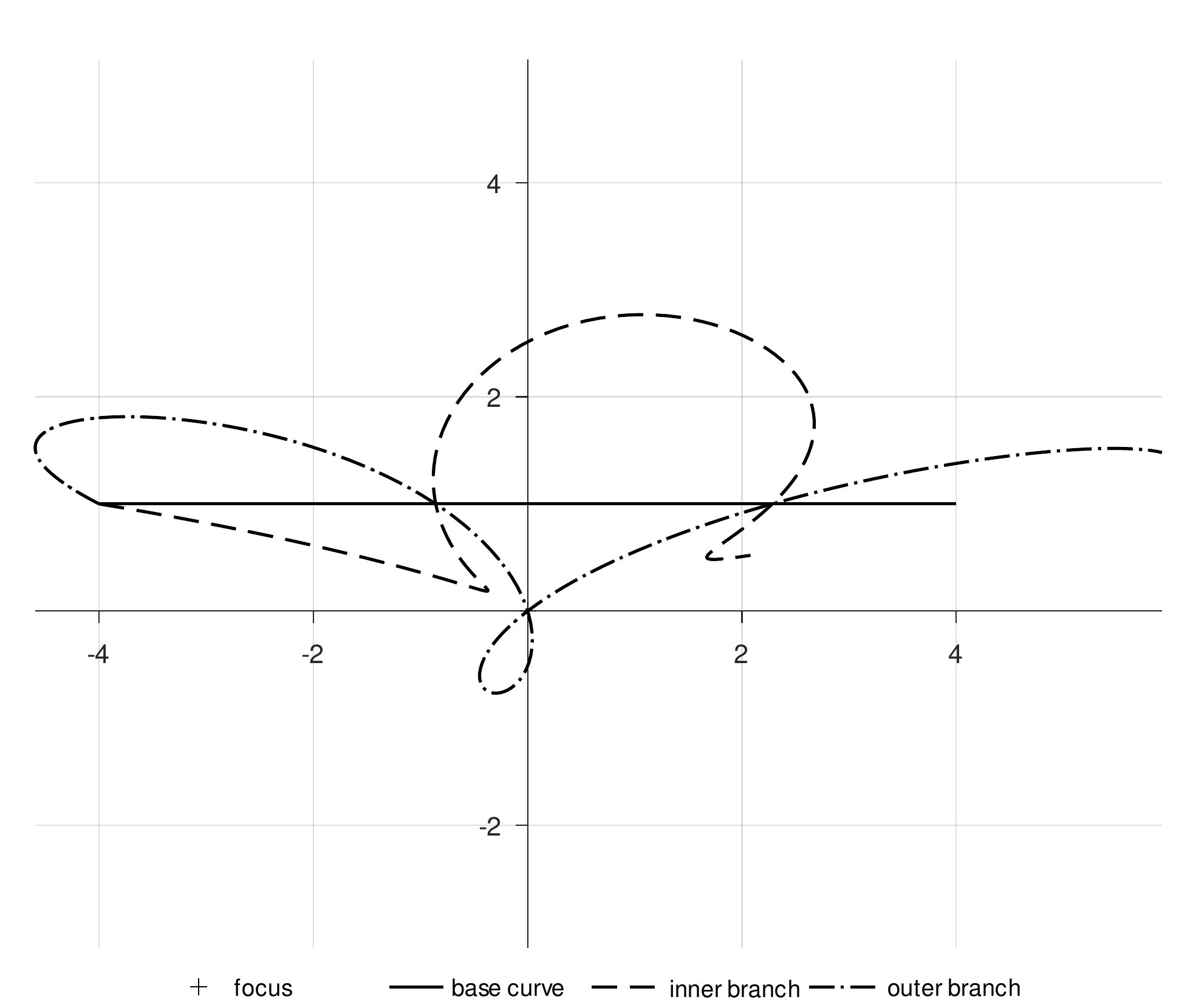}}\quadº
  \subfigure[{\(\mathfrak{C}_{f(l)=\ln{l}}^{O=[0, 0]^\intercal}(\mathcal{C}: y=1)_{N_x=-2  \rightarrow S_x= 2}\)}]{\label{subfig:linConchLogFun} \includegraphics[width=0.45\textwidth]{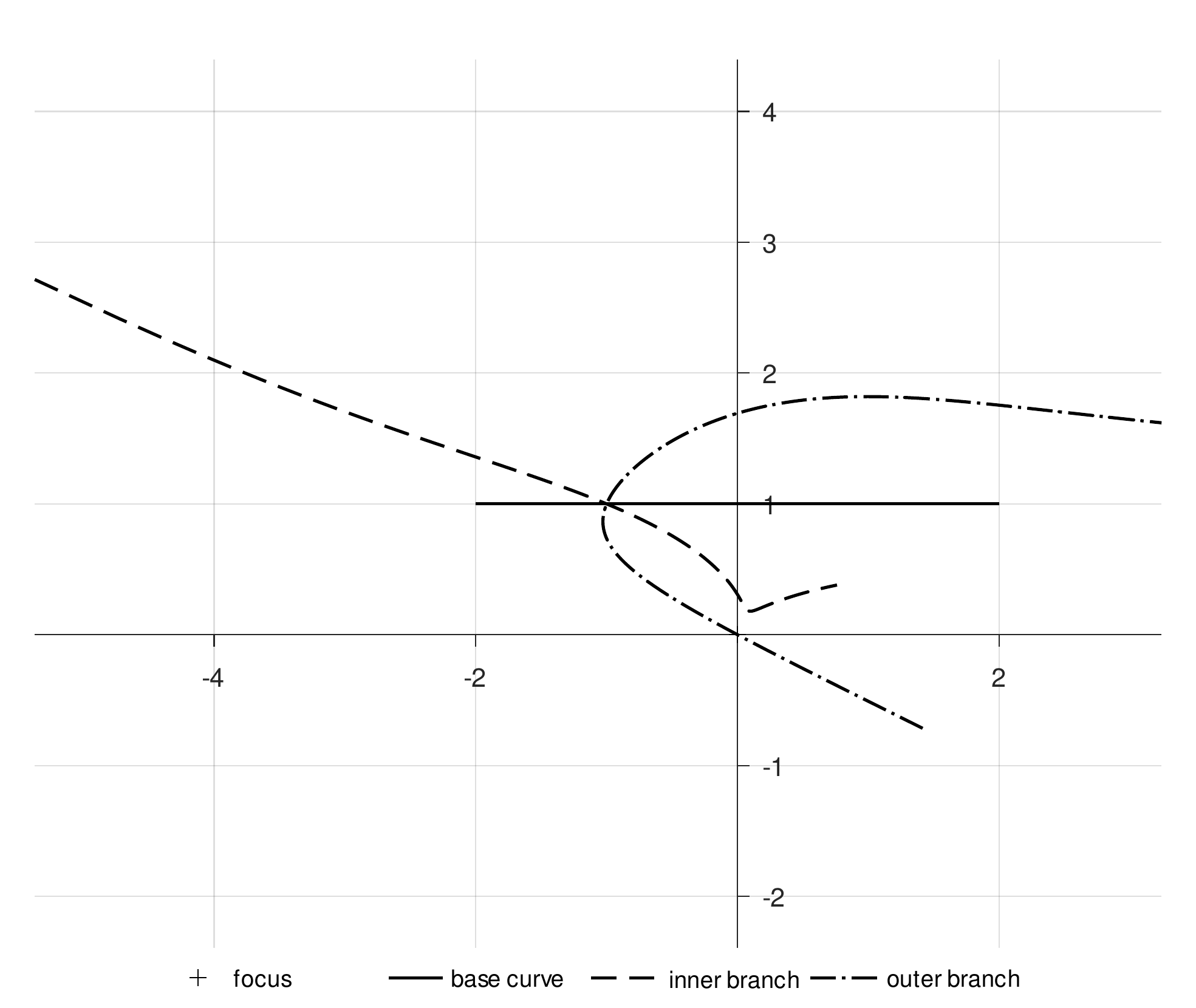}}
\caption{Generalized planar conchoid created when using a line segment as the base curve.}
  \label{fig:GPCwithLine}
\end{figure}

With this implementation, some GPCs were created, as shown in Figure \ref{fig:GPCwithLine}. In the first figure, the \emph{Nicomedes conchoid}, which is a special case of these GPCs, is shown. In this case, the focus is at the origin; \textit{i.e.} \(O=[0, 0]^\intercal\). It has a linear base curve parallel to the x-axis at \(\mathcal{C}: y=1\) from \(N_x=-3\) to \(S_x=3\), and the arc-length function has a constant value; i.e., \(f(l) = 2\).
Other GPCs with the same linear base curves and foci at origin are shown in the figure. The arc length of the GPC shown in the second figure is described with a linear function \(f(l) = l\). By changing the arc-length function, one can obtain a different GPC, as shown in the third figure, where the arc length of the GPC is described with a sine function \(f(l) = \sin{l}\). Similarly, the fourth figure shows a GPC with an arc length described with a logarithmic function \(f(l) = \ln{l}\).

\begin{figure}[!ht]
  \centering
  \subfigure[{\(\mathfrak{C}_{f(l)=\frac{136}{100}}^{O=[0, 0]^\intercal}(\mathcal{C}: c=[0, \frac{113}{100}]^\intercal, r=\frac{80}{100})_{\theta_N=0 \rightarrow \theta_S = 2\pi}\)}]{\label{subfig:circConchConst2Fun} \includegraphics[width=0.45\textwidth]{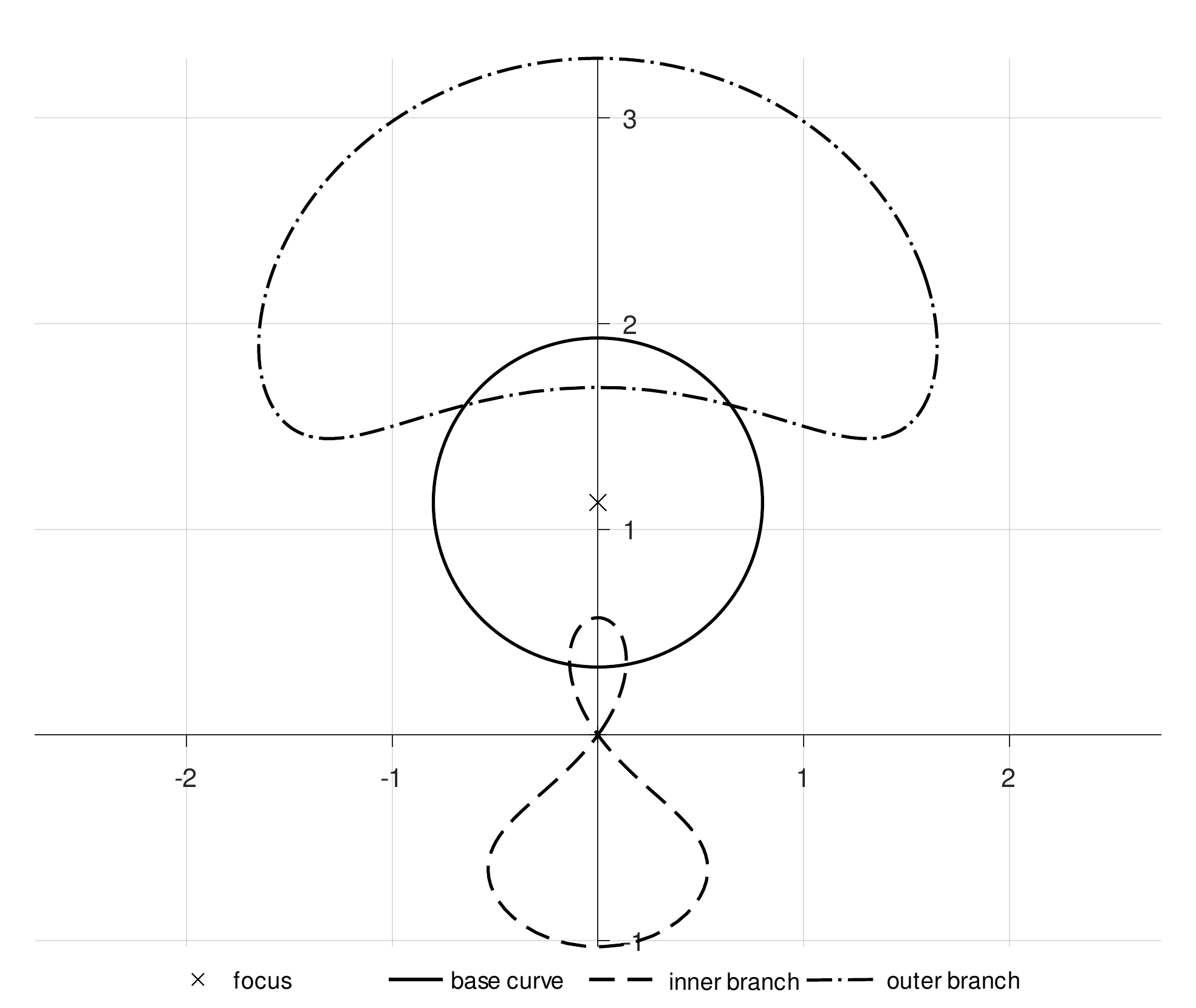}}\quad
  \subfigure[{\(\mathfrak{C}_{f(l)=l}^{O=[0, 0]^\intercal}(\mathcal{C}: c=[0, \frac{7}{2}]^\intercal, r=2)_{\theta_N=0 \rightarrow \theta_S = 2\pi}\)}]{\label{subfig:pcircConchLinFun}
  \includegraphics[width=0.45\textwidth]{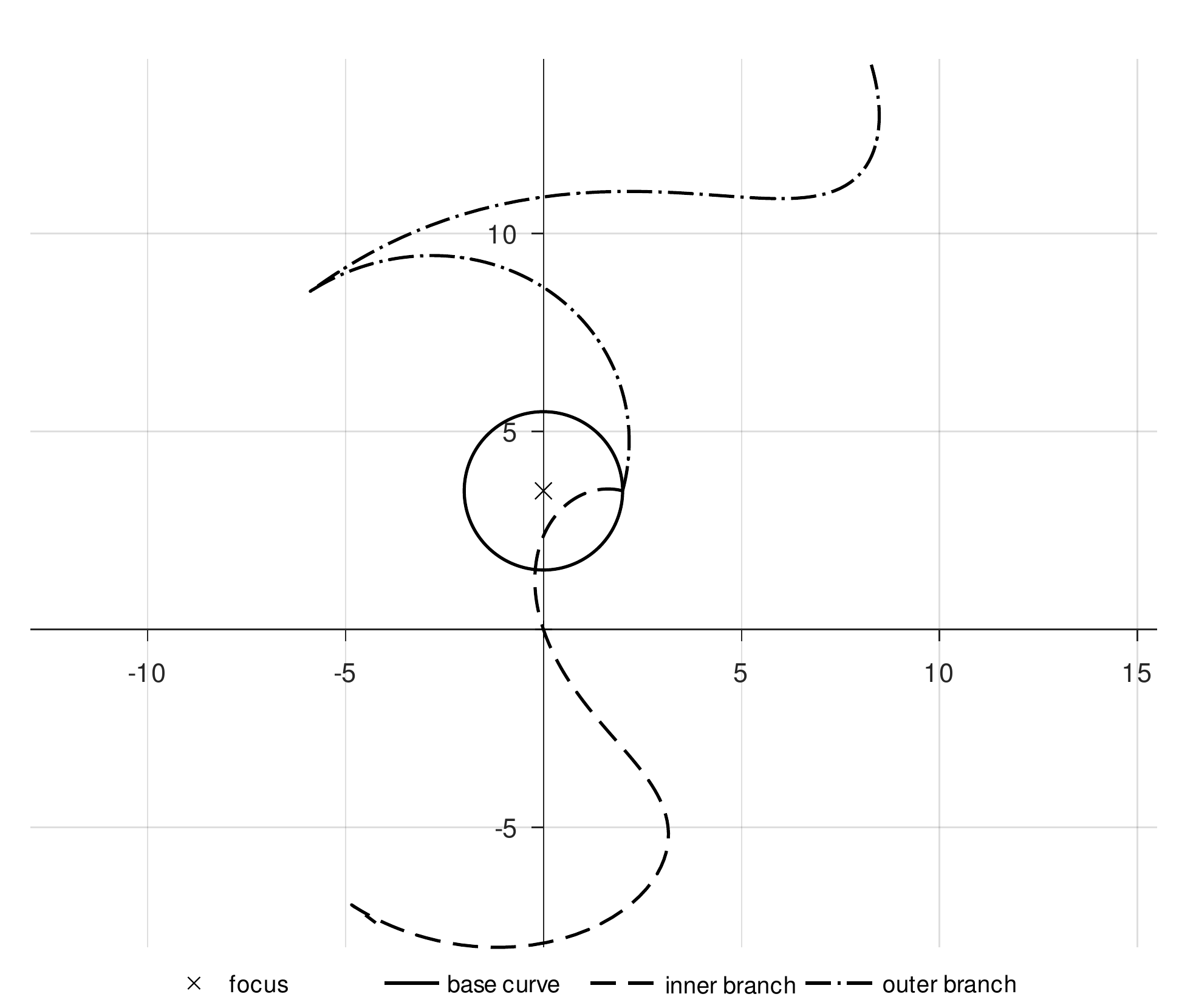}}\\
  \subfigure[{\(\mathfrak{C}_{f(l)=2\sin{l}}^{O=[0, 0]^\intercal}(\mathcal{C}: c=[0, \frac{7}{2}]^\intercal, r=2)_{\theta_N=0 \rightarrow \theta_S = 2\pi}\)}]{\label{subfig:circConchSinFun}
  \includegraphics[width=0.45\textwidth]{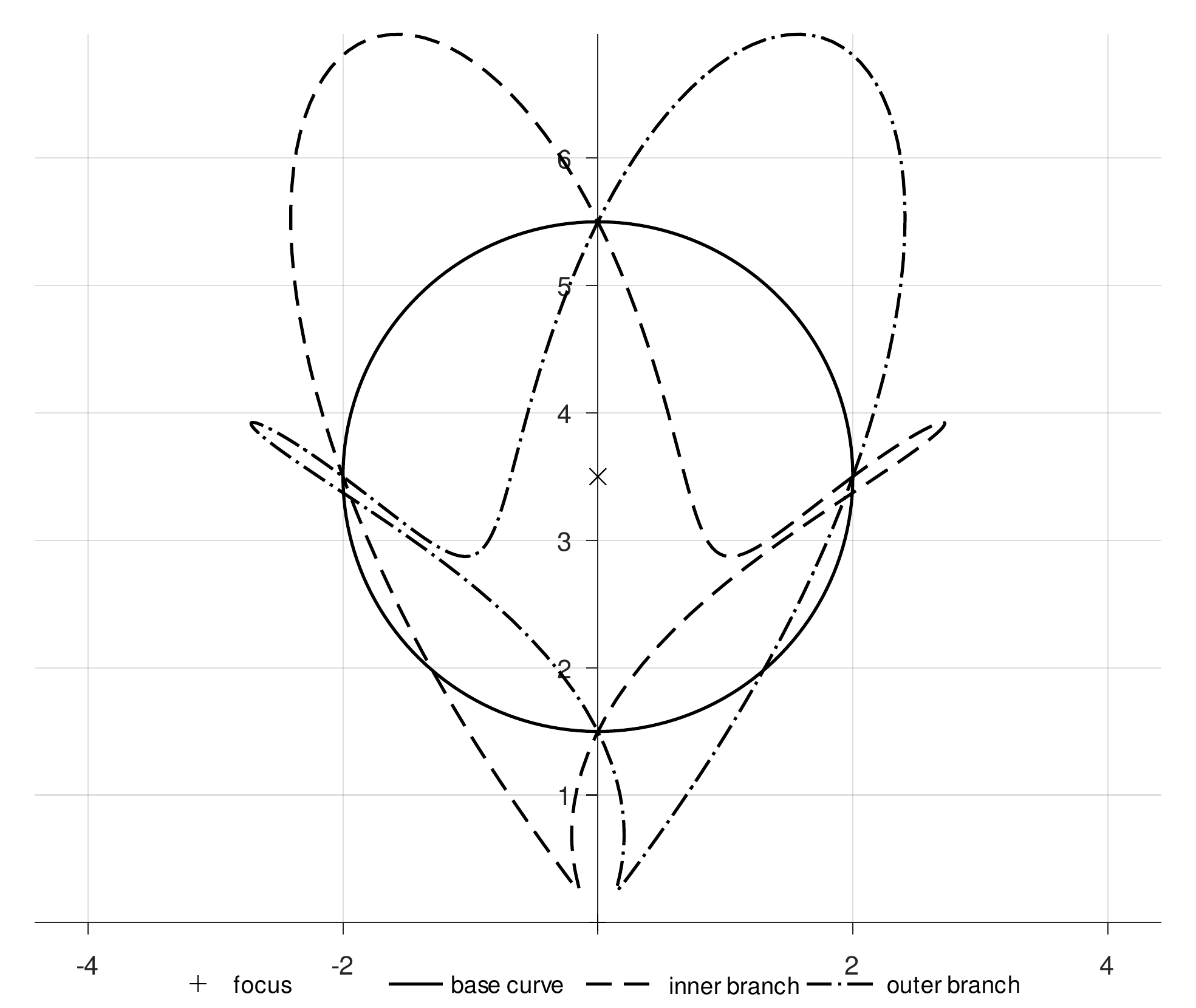}}\quad
  \subfigure[{\(\mathfrak{C}_{f(l)=\ln{l}}^{O=[0, 0]^\intercal}(\mathcal{C}: c=[0, \frac{7}{2}]^\intercal, r=2)_{\theta_N=0 \rightarrow \theta_S = 2\pi}\)}]{\label{subfig:circConchLogFun}
  \includegraphics[width=0.45\textwidth]{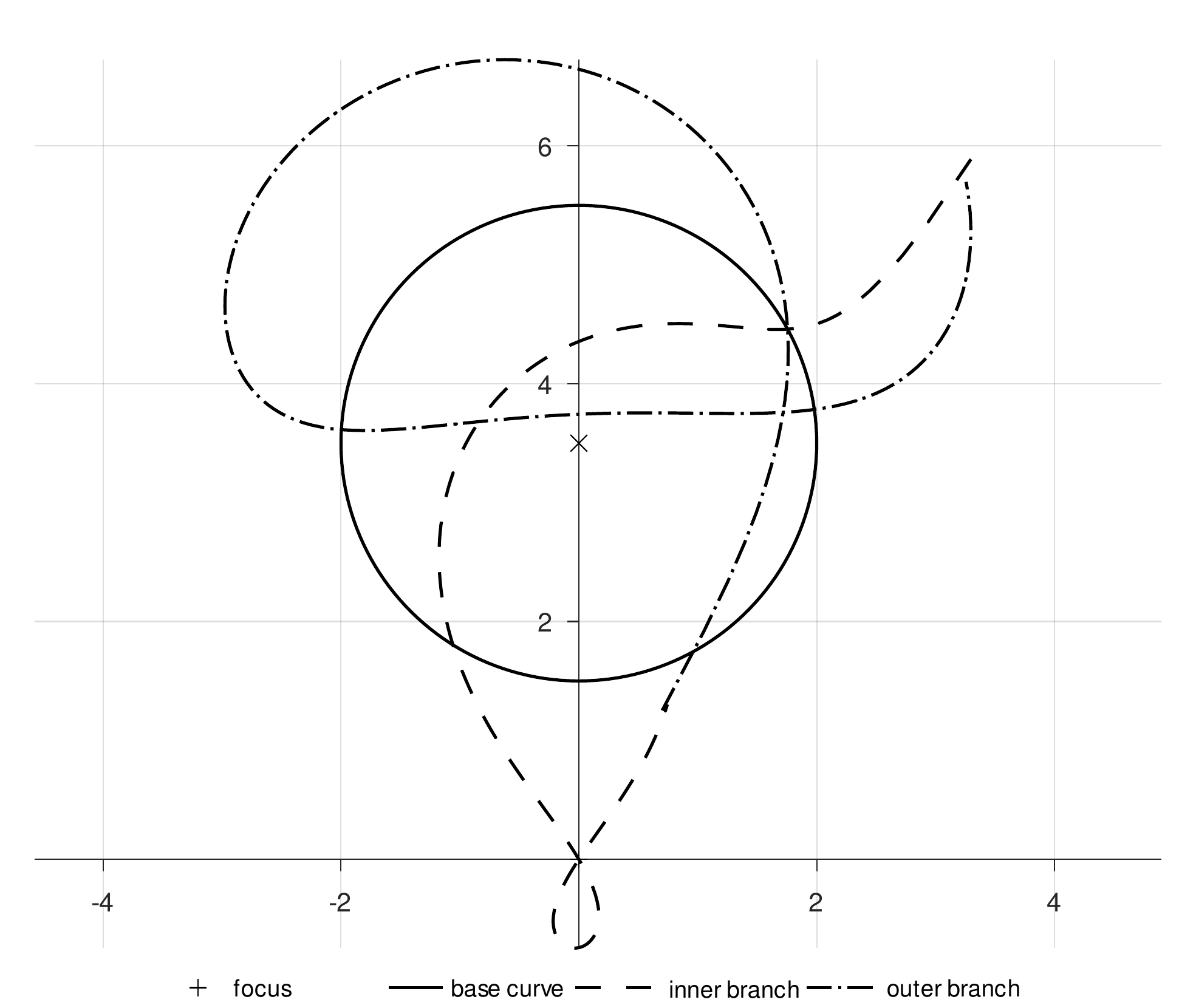}}
  \caption{Generalized planar conchoid created when using an arc circle as a base curve.}
  \label{fig:GPCwithArcCircle}
\end{figure}

Other GPCs were created in Figure \ref{fig:GPCwithArcCircle} by using a circular arc base curve. The first one plots the \emph{Lima\c{c}on of the Pascal conchoid}, which is also a special case of these GPCs. In this particular case, the GPC focus is as the origin of the coordinate system, and the base curve is a circular arc with a centre at \(c\) and a radius of \(r\); \textit{i.e.} \(\mathcal{C}: c=[0, \frac{113}{100}]^\intercal, r=\frac{80}{100}\). The base curve starts when \(\theta_N = 0\) and ends when \(\theta_S = 2\pi\). The arc length function is a constant value of \(f(l) = \frac{136}{100}\). The second
GPC and the subsequent GPCs are generated with the same focus locations. The second GPC, in particular, has a circular arc base curve \(\mathcal{C}: c=[0, \frac{7}{2}]^\intercal, r=2\) from \(\theta_N = 0\) to \(\theta_S = 2\pi\) and an arc-length linear function of \(f(l) = l\). The third GPC is plotted considering a base curve with the properties \(\mathcal{C}: c=[0, \frac{7}{2}]^\intercal, r=2\) and is generated from \(\theta_N = 0\) to \(\theta_S = 2\pi\), and the arc-length function is a trigonometric function; \textit{i.e.} \(f(l) = 2 \sin{l}\). Finally, the last GPC presented here is similar to the last three presented in the figure and is generated with a base curve of \(\mathcal{C}: c=[0, \frac{7}{2}]^\intercal, r=2\) from \(\theta_N = 0\) to \(\theta_S = 2\pi\). Its arc length is bades on a natural logarithm function \(f(l) = \log{l}\).

The computer code scripts used to generate each of the GPCs, including those shown in Figure \ref{fig:GPCwithLine}, are at file \textsf{someLinearGenConchPlotsSCR}, and the script for the GPC presented in Figure \ref{fig:GPCwithArcCircle} is at file \textsf{someCircarcGenConchPlotsSCR}. For example, to generate the third GPC in Figure \ref{fig:GPCwithLine}, the code described in \ref{lst:linearGPC} was used. Similarly, to generate the fourth GPC in Figure \ref{fig:GPCwithArcCircle}, the code described in \ref{lst:arccircleGPC} was used.

\begin{lstlisting}[language=Oz, caption={Code to generate the plot of Figure \ref{fig:GPCwithLine}c}, label={lst:linearGPC}]
clear all
focusColVec = zeros(2,1);
lineFun = @(x) 1;
npts = 180;
abscissaIntval = 4 * [-1, 1];
distFun = @(l) 2 * sin(l);
gcLinPlotStruct = linegenconchoid (focusColVec, lineFun, abscissaIntval,...
distFun, npts);
figure()
hold on
plotgenconchoid (gcLinPlotStruct);
hold off
axis equal
\end{lstlisting}

\begin{lstlisting}[language=Oz, caption={Code to generate the plot of Figure \ref{fig:GPCwithArcCircle}d}, label={lst:arccircleGPC}]
clear all
focusColVec = zeros(2,1);
thetaAngleIntval = [0, 2*pi];
npts = 180;
circStruct = struct("c", [0; 7/2], "r", 2);
distFun = @(l) log(l);
gcCircarcPlotStruct = circarcgenconchoid (focusColVec, circStruct,...
thetaAngleIntval, distFun, npts);
figure()
hold on
plotgenconchoid (gcCircarcPlotStruct);
plot(circStruct.c(1), circStruct.c(2), 'kx');
hold off
axis equal
\end{lstlisting}

Descriptions of the code used in the above scripts are explained in the next section.

\section{The computer code}

The implementation is composed of three functions written in Octave, as can be deduced when reading the listings shown in the preceding paragraph.
The \textsf{linegenconchoid} function generates a GPC with a base curve that is a \emph{line}, while the \textsf{circarcgenconchoid} function generates a GPC with a base curve that is a \emph{circular arc}. Other functions can be created for different base curves. The line and circular-arc functions can be used as templates. Finally, the \textsf{plotgenconchoid} function can be called by any of the above (and newly created) functions to visualize the GPC plot.

The computer code can be downloaded from \cite{SuarezBurgoa2022.code}. In the following subsections, the linear and circular arc functions are described through manual calculations.

To implement the code in MATLAB, only the \textsf{end} statements in every function should be changed. Octave uses particular \textsf{end} statements for each function; for example, the \textsf{if} statement ends with \textsf{endif}. In MATLAB, the ending word for this statement is simply \textsf{end}.

\subsection{The \textsf{linegenconchoid} function}

To show how the \textsf{linegenconchoid} function creates a GPC data structure, we use an example. Our problem is to plot the following GPC:
\[
    \mathfrak{C}_{f(l)= l + \sin{l}}^{O=[2, 1]^\intercal} \left(\mathcal{C}: l = \left[y(x) = \frac{3}{2} + \frac{1}{2} x\right]\right)_{x_N = -3 \rightarrow x_S = 0};
\]
which is read as follows.

A generalized planar conchoid has a focus of
\[
   O = (2, 1),
\]
and \emph{the base curve} is a line given by the function
\[
    y(x) = \frac{3}{2} + \frac{1}{2} x
\]
from point
\[
    N = (-3, y(-3))
\]
to point
\[
    S = (0, y(0)).
\]
The GPC is based the arc-length function
\[
    f(l) = l + \sin{l}.
\]

For convenience, we approximate the GPC with \(m\) discrete points (\(P_i\)) generated in the base curve from point \(N\) to point \(S\); these extreme points are now represented by their corresponding vectors
\[
    N = \boldsymbol{n} = 
   \begin{pmatrix}
       -3\\
       0
   \end{pmatrix}
\]
and
\[
    S = \boldsymbol{s} = 
    \begin{pmatrix}
       0\\
       \frac{3}{2}
   \end{pmatrix}.
\]
Those \(P_i\) points are expressed by vectors \(\boldsymbol{p}_i\)
\[
    \boldsymbol{p}_i = \boldsymbol{n}  + (\boldsymbol{s} - \boldsymbol{n}) k.
\]
The scaling factor \(k\) varies from 0 to 1 and is divided into \(m\) parts
\[
    k = \frac{i}{m},
\]
for \(i = \{0, 2, \ldots, m-1\}\).
The arc length (\(l_i\)) of each line from \(N\) to \(P_i\) is the norm of
\[
    \boldsymbol{l}_i = \boldsymbol{p}_i - \boldsymbol{n},
\]
\textit{id est}
\[
    l_i = |\boldsymbol{l}_i| = \sqrt{\boldsymbol{l}_i \cdot \boldsymbol{l}_i}.
\]

Now, we need to calculate \(d_i = f(l_i)\) to obtain the points that define the inner branch, that is, points \(\boldsymbol{q}_i\). Additionally, with this value, \(d_i\) can generate the points that define the outer branch, that is, points \(\boldsymbol{q}'_i\). We must consider that \(f(l_i)\) is given and that
\[
    d_i = l_i + \sin{l_i}
\]
in this case.
For the current example, Table \ref{tab:ptsBaseCurveArclengthDists} shows the values calculated for the following variables: \(\boldsymbol{p}_i\), \(\boldsymbol{l}_i\), \(|\boldsymbol{l}_i|\), and \(d_i\).

\begin{table}[!ht]
  \label{tab:ptsBaseCurveArclengthDists}
  \caption{Base curve, arc length and distance points.}
  \begin{tabularx}{\textwidth}{X *{6}{p{1.5cm}}}
  \hline
  \multirow{2}{*}{\(k\)} & \multicolumn{2}{c}{\(\boldsymbol{p}_i\)} & \multicolumn{2}{c}{\(\boldsymbol{l}_i\)} & \multirow{2}{*}{\(|\boldsymbol{l}_i|\)} & \multirow{2}{*}{\(d_i\)} \\
  \cline{2-5}
  & \(x\) & \(y\) & \(x\) & \(y\)\\
  \hline
  0.000 & -3.000 & 0.000 & 0.000 & 0.000 & 0.000 & 0.000\\
  0.059 & -2.824 & 0.088 & 0.176 & 0.088 & 0.197 & 0.393\\
  0.118 & -2.647 & 0.176 & 0.353 & 0.176 & 0.395 & 0.779\\
  0.176 & -2.471 & 0.265 & 0.529 & 0.265 & 0.592 & 1.150\\
  0.235 & -2.294 & 0.353 & 0.706 & 0.353 & 0.789 & 1.499\\
  0.294 & -2.118 & 0.441 & 0.882 & 0.441 & 0.987 & 1.821\\
  0.353 & -1.941 & 0.529 & 1.059 & 0.529 & 1.184 & 2.110\\
  0.412 & -1.765 & 0.618 & 1.235 & 0.618 & 1.381 & 2.363\\
  0.471 & -1.588 & 0.706 & 1.412 & 0.706 & 1.578 & 2.578\\
  0.529 & -1.412 & 0.794 & 1.588 & 0.794 & 1.776 & 2.755\\
  0.588 & -1.235 & 0.882 & 1.765 & 0.882 & 1.973 & 2.893\\
  0.647 & -1.059 & 0.971 & 1.941 & 0.971 & 2.170 & 2.996\\
  0.706 & -0.882 & 1.059 & 2.118 & 1.059 & 2.368 & 3.067\\
  0.765 & -0.706 & 1.147 & 2.294 & 1.147 & 2.565 & 3.110\\
  0.824 & -0.529 & 1.235 & 2.471 & 1.235 & 2.762 & 3.133\\
  0.882 & -0.353 & 1.324 & 2.647 & 1.324 & 2.960 & 3.141\\
  0.941 & -0.176 & 1.412 & 2.824 & 1.412 & 3.157 & 3.142\\
  1.000 & 0.000 & 1.500 & 3.000 & 1.500 & 3.354 & 3.143\\
  \hline
  \end{tabularx}
\end{table}

The points \(\boldsymbol{q}_i\) and the locus of the inner branch of the GPC are obtained by a direction unit vector that joins the focus (point \(O\)) with point \(P_i\) and the distance \(d_i\) to the inside from \(O\).
The unitary vector for each \(\boldsymbol{p}_i\) from \(O\) is
\begin{eqnarray}
    \boldsymbol{u}_o &=& \frac{\boldsymbol{p}_i - \boldsymbol{o}}{|\boldsymbol{p}_i - \boldsymbol{o}|}, \\
    &=&
    \frac{\boldsymbol{p}_i - \boldsymbol{o}}{\sqrt{(\boldsymbol{p}_i - \boldsymbol{o}) \cdot (\boldsymbol{p}_i - \boldsymbol{o})}}.
\end{eqnarray}
Table \ref{tab:vectsPOandUniVect} shows the intermediate variable values needed to obtain the unitary vector \(\boldsymbol{u}_o\).

\begin{table}[!ht]
  \label{tab:vectsPOandUniVect}
  \caption{Vector \((\boldsymbol{p}_i - \boldsymbol{o})\) and its unit vector \(\boldsymbol{u}_o\).}
  \begin{tabularx}{\textwidth}{X *{5}{p{1.8cm}}}
  \hline
  \multirow{2}{*}{\(k\)} & \multicolumn{2}{c}{\((\boldsymbol{p}_i - \boldsymbol{o})\)} & \multirow{2}{*}{\(|\boldsymbol{p}_i - \boldsymbol{o}|\)} & \multicolumn{2}{c}{\(\boldsymbol{u}_o\)} \\
  \cline{2-3}\cline{5-6}
  & \(x\) & \(y\) & & \(x\) & \(y\)\\
  \hline
  0.000 & -3.000 & 0.000 & 3.000 & -1.000 & 0.000\\
  0.059 & -2.824 & 0.088 & 2.825 & -1.000 & 0.031\\
  0.118 & -2.647 & 0.176 & 2.653 & -0.998 & 0.067\\
  0.176 & -2.471 & 0.265 & 2.485 & -0.994 & 0.107\\
  0.235 & -2.294 & 0.353 & 2.321 & -0.988 & 0.152\\
  0.294 & -2.118 & 0.441 & 2.163 & -0.979 & 0.204\\
  0.353 & -1.941 & 0.529 & 2.012 & -0.965 & 0.263\\
  0.412 & -1.765 & 0.618 & 1.870 & -0.944 & 0.330\\
  0.471 & -1.588 & 0.706 & 1.738 & -0.914 & 0.406\\
  0.529 & -1.412 & 0.794 & 1.620 & -0.872 & 0.490\\
  0.588 & -1.235 & 0.882 & 1.518 & -0.814 & 0.581\\
  0.647 & -1.059 & 0.971 & 1.436 & -0.737 & 0.676\\
  0.706 & -0.882 & 1.059 & 1.378 & -0.640 & 0.768\\
  0.765 & -0.706 & 1.147 & 1.347 & -0.524 & 0.852\\
  0.824 & -0.529 & 1.235 & 1.344 & -0.394 & 0.919\\
  0.882 & -0.353 & 1.324 & 1.370 & -0.258 & 0.966\\
  0.941 & -0.176 & 1.412 & 1.423 & -0.124 & 0.992\\
  1.000 & 0.000 & 1.500 & 1.500 & 0.000 & 1.000\\
  \hline
  \end{tabularx}
\end{table}

With the above calculated variables, we can calculate the points that define the locus of the inner branch of the GPC, which is obtained with the equation
\[
    \boldsymbol{q}_i = \boldsymbol{p}_i - d_i \boldsymbol{u}_o.
\]
Similarly, for the case of the points \(\boldsymbol{q'}_i\) that define the locus of the outer branch of the GPC, the equation is
\[
    \boldsymbol{q'}_i = \boldsymbol{p}_i + d_i \boldsymbol{u}_o.
\]
These calculated coordinates are shown in Table \ref{tab:ptsInnerOuterBranches}.

\begin{table}[!ht]
  \label{tab:ptsInnerOuterBranches}
  \caption{Points on the inner and outer branches.}
  \begin{tabularx}{\textwidth}{X *{4}{p{2.2cm}}}
  \hline
  \multirow{2}{*}{\(k\)} & \multicolumn{2}{c}{\(\boldsymbol{q}_i\)} & \multicolumn{2}{c}{\(\boldsymbol{q}'_i\)} \\
  \cline{2-5}
  & \(x\) & \(y\) & \(x\) & \(y\)\\
  \hline
  0.000 & -3.000 & 0.000 & -3.000 & 0.000\\
  0.059 & -2.430 & 0.076 & -3.217 & 0.101\\
  0.118 & -1.870 & 0.125 & -3.424 & 0.228\\
  0.176 & -1.327 & 0.142 & -3.614 & 0.387\\
  0.235 & -0.813 & 0.125 & -3.776 & 0.581\\
  0.294 & -0.335 & 0.070 & -3.900 & 0.812\\
  0.353 & 0.094 & -0.026 & -3.977 & 1.085\\
  0.412 & 0.466 & -0.163 & -3.995 & 1.398\\
  0.471 & 0.768 & -0.341 & -3.944 & 1.753\\
  0.529 & 0.989 & -0.556 & -3.813 & 2.145\\
  0.588 & 1.119 & -0.799 & -3.590 & 2.564\\
  0.647 & 1.150 & -1.054 & -3.267 & 2.995\\
  0.706 & 1.081 & -1.297 & -2.846 & 3.415\\
  0.765 & 0.924 & -1.502 & -2.336 & 3.796\\
  0.824 & 0.705 & -1.644 & -1.763 & 4.115\\
  0.882 & 0.456 & -1.711 & -1.162 & 4.358\\
  0.941 & 0.213 & -1.706 & -0.566 & 4.529\\
  1.000 & 0.000 & -1.643 & 0.000 & 4.643\\
  \hline
  \end{tabularx}
\end{table}

Finally, the GPC we are looking for is approximated by the set of discrete points \(\boldsymbol{q}_i\) and \(\boldsymbol{q'}_i\), i.e.,
\[
    \mathfrak{C} = \{\boldsymbol{q}_i, \boldsymbol{q'}_i\}.
\]
The GPC that is calculated the step-by-step using intermediate variables to finally obtain the coordinates is shown in Table \ref{tab:ptsInnerOuterBranches}.

\subsection{The \textsf{circarcgenconchoid} function}

Similar to the preceding case, to show how the \textsf{circarcgenconchoid} function creates the GPC data structure, we make use of another example that is translated for solving the problem of plotting the GPC; it is written as
\[
    \mathfrak{C}_{f(l)=l + \frac{1}{l}}^{O=[0, 0]^\intercal}(\mathcal{C}: c=[5, 10]^\intercal, r=6)_{\theta_N=0 \rightarrow \theta_S = \frac{9}{8}\pi}
\]
which is read as follows.

A generalized planar conchoid has a focus of
\[
   O = (0, 0),
\]
and a \emph{base curve} that is a circular arc with a centre \(C\) and radius \(r\), which are equal to
\[
    C = (5, 10)
\]
\[
    r = 6,
\]
that starts from an angle of
\[
    \theta_{N} = 0
\]
and ends at an angle of
\[
    \theta_{S} = \frac{9}{8} \pi
\]
is based on the arc-length function
\[
    f(l) = l + \frac{1}{l}.
\]
Here, we show how to approximate this GPC with \(m=180\) \(P_i\) points generated on the base curve from \(N\) to \(S\).

The circle to which the circular arc belongs is defined by its centre (now as a vector \(\boldsymbol{c}\)) and radius \(r\). As defined by a polar equation, the circle's function is
\[
    \begin{pmatrix}
        x \\ y
    \end{pmatrix}
    =
    \boldsymbol{c} + r
    \begin{pmatrix}
        \cos{\theta} \\ \sin{\theta}
    \end{pmatrix}.
\]

Then, points \(N\) and \(S\), which are represented by vectors in this case, are
\begin{eqnarray}
    N = \boldsymbol{n} &=& 
   \begin{pmatrix}
       5\\
       10
   \end{pmatrix}
   + 6
   \begin{pmatrix}
       \cos{0}\\
       \sin{0}
   \end{pmatrix} \\
   &=&
   \begin{pmatrix}
       11\\
       10
   \end{pmatrix}
\end{eqnarray}
and
\begin{eqnarray}
    S = \boldsymbol{n} &=& 
   \begin{pmatrix}
       5\\
       10
   \end{pmatrix}
   + 6
   \begin{pmatrix}
       \cos{\left(\frac{9}{8}\pi\right)}\\
       \sin{\left(\frac{9}{8}\pi\right)}
   \end{pmatrix} \\
   &=&
   \begin{pmatrix}
       5 - 3 \sqrt{2 + \sqrt{2}}\\
       10 - 3 \sqrt{2 - \sqrt{2}}
   \end{pmatrix} \\
   &\approx&
   \begin{pmatrix}
       -0.543\\
        7.704
   \end{pmatrix}.   
\end{eqnarray}

The arc length between these two points is
\[
    L = r (\theta_S - \theta_N).
\]

Every point at the base curve \(\boldsymbol{p}_i\) is distributed on the circle between \(N\) and \(S\); then,
\[
    \boldsymbol{p}_i = \boldsymbol{c} + r
    \begin{pmatrix}
        \cos{\theta_i} \\ \sin{\theta_i}
    \end{pmatrix};
\]
where
\[
    \theta_i = \theta_N + \frac{L}{r} k
\]
and the scaling factor \(k\) varies from 0 to 1 and is divided into \(m\) parts; \textit{i.e.}
\[
    k = \frac{i}{m},
\]
for \(i = \{0, 2, \ldots, m-1\}\).
The arc lengths from \(N\) to every \(\boldsymbol{p}_i\) are
\[
    l_i = L\; k.
\]

The following examples are similar to those described in the preceding subsection when using the linear base curve; because now, we have an expression for \(l_i\), which can be passed through the function \(f(l_i)\) in this example.
\[
    f(l) = l + \frac{1}{l};
\]
Therefore,
\[
   d_i = l_i + \frac{1}{l_i}.
\]
This finally gives that
\[
    \boldsymbol{q}_i = \boldsymbol{p}_i - d_i \boldsymbol{u}_o.
\]
and
\[
    \boldsymbol{q'}_i = \boldsymbol{p}_i + d_i \boldsymbol{u}_o.
\]
The same expression is used for \(\boldsymbol{u}_o\); \textit{i.e.}
\[
    \boldsymbol{u}_o = 
    \frac{\boldsymbol{p}_i - \boldsymbol{o}}{\sqrt{(\boldsymbol{p}_i - \boldsymbol{o}) \cdot (\boldsymbol{p}_i - \boldsymbol{o})}}.
\]

The GPC we are looking for is approximated by the sets of discrete points \(\boldsymbol{q}_i\) and \(\boldsymbol{q'}_i\)
\[
    \mathfrak{C} = \{\boldsymbol{q}_i, \boldsymbol{q'}_i\}.
\]
Indeed, this GPC is the same as the GPC plotted in Figure \ref{fig:partsOfGCdefinition}. This time, it was necessary to further discretize the curve to find a good approximation (180 points were used). For that reason, a step-by-step calculation table is not presented, but the script shown in \ref{lst:gpcCircArcLst3} can be used to visualize this result.

\begin{lstlisting}[language=Oz, caption={Code to generate Figure \ref{fig:partsOfGCdefinition}b}, label={lst:gpcCircArcLst3}]
clear all
## comment
focusColVec = zeros(2,1);
thetaAngleIntval = [0*pi, 9/8*pi];
npts = 18*10;

circStruct = struct("c", [5; 10], "r", 6);
distFun = @(l) 1./l + l;
gcCircarcPlotStruct = circarcgenconchoid (focusColVec, circStruct,...
thetaAngleIntval, distFun, npts);

figure()
hold on
plotgenconchoid (gcCircarcPlotStruct);
plot(circStruct.c(1), circStruct.c(2), 'kx');
hold off
axis equal
\end{lstlisting}

\section{Final remark}

The user can create many GPCs depending on the different focus places, base curves, arc-length functions and intervals.





\renewcommand\refname{References}
\bibliographystyle{chicagoa}
\bibliography{references.bib}

\begin{thebibliography}{}

\bibitem[\protect\citeauthoryear{{Su\'arez-Burgoa}}{{Su\'arez-Burgoa}}{2022}]{SuarezBurgoa2022.code}
{Su\'arez-Burgoa}, L. (2022, Apr).
\newblock gen\uppercase{P}lanar\uppercase{C}onchoid: Computer code in
  \uppercase{MATLAB}/\uppercase{O}ctave for the calculation of a generalized
  planar conchoid for two general base-curves, linear-segment and circular-arc.
\newblock \url{https://github.com/losuarezburgoa/genPlanarConchoid}.

\end{thebibliography}



\end{document}